\input amstex
\documentstyle{amsppt}
%
\catcode`@=11
\redefine\output@{%
  \def\break{\penalty-\@M}\let\par\endgraf
  \ifodd\pageno\global\hoffset=105pt\else\global\hoffset=8pt\fi  
  \shipout\vbox{%
    \ifplain@
      \let\makeheadline\relax \let\makefootline\relax
    \else
      \iffirstpage@ \global\firstpage@false
        \let\rightheadline\frheadline
        \let\leftheadline\flheadline
      \else
        \ifrunheads@ 
        \else \let\makeheadline\relax
        \fi
      \fi
    \fi
    \makeheadline \pagebody \makefootline}%
  \advancepageno \ifnum\outputpenalty>-\@MM\else\dosupereject\fi
}
\def\Beta{\mathchar"0\hexnumber@\rmfam 42}
\catcode`\@=\active
\nopagenumbers
\chardef\textvolna='176

\font\tencyr=wncyr10
\font\tenbfit=cmbxti10
\chardef\bigalpha='013
\def\negskp{\hskip -2pt}

\def\compos{\,\raise 1pt\hbox{$\sssize\circ$} \,}

\chardef\degree="5E
\accentedsymbol\tildepsi{\kern 3.2pt\tilde{\kern -3.2pt\psi}}
\accentedsymbol\ssizetildepsi{\ssize\kern 2.5pt\tilde{\kern -2.5pt\psi}}
\def\Pover#1{\overset{\kern 1.5pt #1}\to P}
\def\Qover#1{\overset{\kern 1.5pt #1}\to Q}
\def\Rover#1{\overset{\kern 1.5pt #1}\to R}
\def\Sover#1{\overset{\kern 1.5pt #1}\to S}
\def\blue#1{#1}

\catcode`#=11\def\diez{#}\catcode`#=6
\catcode`&=11\catcode`&=4
\catcode`_=11\catcode`_=8
\catcode`~=11\catcode`~=\active
\def\mycite#1{\cite{\blue{#1}}\immediate\special{ps:
     ShrHPSdict begin /ShrBORDERthickness 0 def}}
\def\myciterange#1#2#3#4{\cite{\blue{#2#3#4}}\immediate\special{ps:
     ShrHPSdict begin /ShrBORDERthickness 0 def}}
\def\mytag#1{%
    \tag#1}
\def\mythetag#1{\thetag{\blue{#1}}\immediate\special{ps:
     ShrHPSdict begin /ShrBORDERthickness 0 def}}
\def\myrefno#1{\no#1}
\def\myhref#1#2{\blue{#2}\immediate\special{ps:
     ShrHPSdict begin /ShrBORDERthickness 0 def}}
\def\myEarXivlink{\myhref{http://arXiv.org}{http:/\negskp/arXiv.org}}

\def\mytheorem#1{\csname proclaim\endcsname{Theorem #1}}
\def\mytheoremwithtitle#1#2{\csname proclaim\endcsname{Theorem #1#2}}
\def\mythetheorem#1{\blue{#1}\immediate\special{ps:
     ShrHPSdict begin /ShrBORDERthickness 0 def}}
\def\mylemma#1{\csname proclaim\endcsname{Lemma #1}}
\def\mylemmawithtitle#1#2{\csname proclaim\endcsname{Lemma #1#2}}
\def\mythelemma#1{\blue{#1}\immediate\special{ps:
     ShrHPSdict begin /ShrBORDERthickness 0 def}}
\def\mycorollary#1{\csname proclaim\endcsname{Corollary #1}}

\def\mydefinition#1{\definition{Definition #1}}
\def\mythedefinition#1{\blue{#1}\immediate\special{ps:
     ShrHPSdict begin /ShrBORDERthickness 0 def}}
\def\myconjecture#1{\csname proclaim\endcsname{Conjecture #1}}
\def\myconjecturewithtitle#1#2{\csname proclaim\endcsname{Conjecture #1#2}}

\def\myproblem#1{\csname proclaim\endcsname{Problem #1}}
\def\myproblemwithtitle#1#2{\csname proclaim\endcsname{Problem #1#2}}


\pagewidth{360pt}
\pageheight{606pt}
\topmatter
\title
Comparison of two classifications of a class of ODE's
in the case of general position.
\endtitle
\rightheadtext{Comparison of two classifications.}
\author
Ruslan Sharipov
\endauthor
\address Bashkir State University, 32 Zaki Validi street, 450074 Ufa, Russia
\endaddress
\email\myhref{mailto:r-sharipov\@mail.ru}{r-sharipov\@mail.ru}
\endemail
\abstract
     Two classifications of second order ODE's cubic with respect to the 
first order derivative are compared in the case of general position, which is
common for both classifications. The correspondence of vectorial, pseudovectorial,
scalar, and pseudoscalar invariants is established. 
\endabstract
\endtopmatter
\loadbold
\TagsOnRight
\document

\head
1. Introduction.
\endhead
     The main object for the present research is the following class of ordinary 
differential equations cubic with respect to the first order derivative:
$$
\hskip -2em
y''=P(x,y)+3\,Q(x,y)\,y'+3\,R(x,y)\,(y')^2+S(x,y)\,(y')^3.
\mytag{1.1}
$$
Differential equations of the form \mythetag{1.1} attracted the attention since 
the epoch of classical papers (see \mycite{1} and \mycite{2}). Nowadays they 
are studied in a large variety of papers so that one 
cannot cite all of them without risks to unintentionally miss some names. With my 
apologies I proceed to the papers \mycite{3} and \mycite{4}, where the equations 
\mythetag{1.1} were classified. They were subdivided into {\bf nine subclasses} 
invariant with respect to transformations of the form 
$$
\hskip -2em
\cases
\tilde x=\tilde x(x,y),\\
\tilde y=\tilde y(x,y).
\endcases
\mytag{1.2}
$$
This subdivision is based on scalar invariants of the equations \mythetag{1.1}
(see definition below) and on their symmetry groups, which are nontrivial
in some cases.\par
     Recently Yu\.~Yu\.~Bagderina in \mycite{5} presented another classification
of the equations \mythetag{1.1} again subdividing them into {\bf nine subclasses} 
invariant with respect to transformations of the form \mythetag{1.2} (see Theorem 2
in \mycite{5}). Her approach is based on Sophus Lie's method of infinitesimal 
transformations adapted to equations of the form \mythetag{1.1} by N\.~H\.~Ibragimov 
in \mycite{6}.\par 
    Unfortunately in \mycite{5} Yu\.~Yu\.~Bagderina does not mention the previously
existing classification from \mycite{3} and \mycite{4}. Though she cites the paper 
\mycite{4} in \mycite{5}, she references this paper only as a source of invariants 
and for criticism of its method.\par
    In the present paper I start a research intended to examine the results 
of paper \mycite{5} and compare them with the prior results from 
\myciterange{3}{3}{, }{4} and \mycite{7}. 
\head
2. Some notations and definitions.
\endhead
    Transformations \mythetag{1.2} are called point transformations. They are
assumed to be locally invertible. For their inverse transformations we write
$$
\hskip -2em
\cases
x=\tilde x(\tilde x,\tilde y),\\
y=\tilde y(\tilde x,\tilde y).
\endcases
\mytag{2.1}
$$
Then, according to \myciterange{3}{3}{, }{4} and \mycite{7}, here we introduce 
the following notations for partial derivatives of the functions in 
\mythetag{1.2} and \mythetag{2.1}:
$$
\xalignat 4
&\hskip -2em
\tilde x_{\sssize 1.0}=\frac{\partial\tilde x}{\partial x},
&&\tilde x_{\sssize 0.1}=\frac{\partial\tilde x}{\partial y},
&&\tilde y_{\sssize 0.1}=\frac{\partial\tilde y}{\partial x},
&&\tilde y_{\sssize 0.1}=\frac{\partial\tilde y}{\partial y},
\qquad
\mytag{2.2}\\
\vspace{2ex}
&\hskip -2em
x_{\sssize 1.0}=\frac{\partial x}{\partial\tilde x},
&&x_{\sssize 0.1}=\frac{\partial x}{\partial\tilde y},
&&y_{\sssize 0.1}=\frac{\partial y}{\partial\tilde x},
&&y_{\sssize 0.1}=\frac{\partial y}{\partial\tilde y}.
\qquad
\mytag{2.3}
\endxalignat
$$
In terms of the notations \mythetag{2.2} and \mythetag{2.3} the Jacoby 
matrices of the point transformations \mythetag{1.2} and \mythetag{2.1} 
are written as follows 
$$
\xalignat 2
&\hskip -2em
S=\Vmatrix
x_{\sssize 1.0} &x_{\sssize 0.1}\\
\vspace{1ex}
y_{\sssize 1.0} &y_{\sssize 0.1}
\endVmatrix,
&&T=\Vmatrix
\tilde x_{\sssize 1.0} &\tilde x_{\sssize 0.1}\\
\vspace{1ex}
\tilde y_{\sssize 1.0} &\tilde y_{\sssize 0.1}
\endVmatrix.
\mytag{2.4}
\endxalignat
$$
In differential geometry the matrices \mythetag{2.1} are called the direct 
and inverse transition matrices (see \mycite{8}). As for the transformations
\mythetag{1.2} and \mythetag{2.1}, geometrically they are interpreted as
changes of local curvilinear coordinates on the plane $\Bbb R^2$ or on some
two-dimensional manifold.\par
     Tensorial and pseudotensorial fields in local coordinates are presented
by their components forming arrays of functions whose arguments are $x,\,y$
or $\tilde x,\,\tilde y$ respectively. In the coordinate form they can be 
defined as follows. 
\mydefinition{2.1} A pseudotensorial field of the type $(r,s)$ and weight $m$ 
is an array of quantities $F^{i_1\ldots\,i_r}_{j_1\ldots\,j_s}$ which under 
the change of coordinates \mythetag{1.2} transforms as 
$$
\hskip -2em
F^{i_1\ldots\,i_r}_{j_1\ldots\,j_s}=
(\det T)^m\sum\Sb p_1\ldots p_r\\ q_1\ldots q_s\endSb
S^{i_1}_{p_1}\ldots\,S^{i_r}_{p_r}\,\,
T^{q_1}_{j_1}\ldots\,T^{q_s}_{j_s}\,\,
\tilde F^{p_1\ldots\,p_r}_{q_1\ldots\,q_s}.
\mytag{2.5}
$$
Tensorial fields are those pseudotensorial fields whose weight is zero, i\.\,e\.
$m=0$ in \mythetag{2.5} for a tensorial field.
\enddefinition
\mydefinition{2.2} Vectorial and pseudovectorial fields are those fields in
Definition~\mythedefinition{2.1} whose type is $(1,0)$. Covectorial and 
pseudocovectorial fields are those fields in Definition~\mythedefinition{2.1}
whose type is $(0,1)$. Scalar and pseudoscalar fields are those fields 
whose type is $(0,0)$.
\enddefinition
\mydefinition{2.3} Tensorial and pseudotensorial fields whose components are 
expressed through $y'$, through the coefficients $P$, $Q$, $R$, $S$ of the equation
\mythetag{1.1}, and through their partial derivatives are called tensorial and 
pseudotensorial invariants of this equation respectively.
\enddefinition
     Invariants of the equation \mythetag{1.1} are subdivided into absolute
and relative ones. \pagebreak Tensorial invariants are absolute invariants, while 
pseudotensorial invariants are relative ones. This definition of relative 
invariants is close to that of \mycite{9}. In some papers (e\.\,g\. in \mycite{5}) 
the term absolute  invariant is applied to scalar invariants only. The definition 
of relative invariants in \mycite{5} is rather loose. Therefore the usage of this 
term in \mycite{5} is different from the above definition.\par
\head
3. Comparison of some relative invariants.
\endhead	
     In \mycite{5} Yu\.~Yu\. Bagderina introduces a long list of special 
notations. She call them relative invariants of various orders from one
to six. The first order expressions are given by the formulas \thetag{2.1} 
in \mycite{5}. Here are they:
$$
\hskip -2em
\aligned
&\alpha^{\sssize\text{Bgd}}_0=Q_{\sssize 1.0}-P_{\sssize 0.1}+2\,P\,R-2\,Q^2,\\
&\alpha^{\sssize\text{Bgd}}_1=R_{\sssize 1.0}-Q_{\sssize 0.1}+P\,S-Q\,R,\\
&\alpha^{\sssize\text{Bgd}}_2=S_{\sssize 1.0}-R_{\sssize 0.1}+2\,Q\,S-2\,R^2.
\endaligned
\mytag{3.1}
$$
Some quantities by Yu\.~Yu\. Bagderina from \mycite{5} share the same symbol
with absolutely different quantities from \myciterange{3}{3}{, }{4} and 
\mycite{7}. Therefore in \mythetag{3.1} and in other formulas I use the upper 
mark {\tencyr<}Bgd{\tencyr>} in order to distinguish Bagderina's quantities 
from those of me and Vera~V.~Kartak (Dmitrieva).\par
     The second order expressions are given by the formulas \thetag{2.2} 
in \mycite{5}: 
$$
\hskip -2em
\aligned
&\beta^{\,\sssize\text{Bgd}}_1=\partial_x\alpha^{\sssize\text{Bgd}}_1
-\partial_y\alpha^{\sssize\text{Bgd}}_0+R\,\alpha^{\sssize\text{Bgd}}_0
-2\,Q\,\alpha^{\sssize\text{Bgd}}_1+P\,\alpha^{\sssize\text{Bgd}}_2,\\
&\beta^{\,\sssize\text{Bgd}}_2=\partial_x\alpha^{\sssize\text{Bgd}}_2
-\partial_y\alpha^{\sssize\text{Bgd}}_1+S\,\alpha^{\sssize\text{Bgd}}_0
-2\,R\,\alpha^{\sssize\text{Bgd}}_1+Q\,\alpha^{\sssize\text{Bgd}}_2.
\endaligned
\mytag{3.2}
$$
The third order expressions are given by the formulas \thetag{2.3} 
in \mycite{5}:  
$$
\hskip -2em
\aligned
&\gamma^{\sssize\text{Bgd}}_{10}=\partial_x\beta^{\,\sssize\text{Bgd}}_1
-Q\,\beta^{\,\sssize\text{Bgd}}_1+P\,\beta^{\,\sssize\text{Bgd}}_2,\\
&\gamma^{\sssize\text{Bgd}}_{11}=\partial_x\beta^{\,\sssize\text{Bgd}}_2
-R\,\beta^{\,\sssize\text{Bgd}}_1+Q\,\beta^{\,\sssize\text{Bgd}}_2,\\
&\gamma^{\sssize\text{Bgd}}_{20}=\partial_y\beta^{\,\sssize\text{Bgd}}_1
-R\,\beta^{\,\sssize\text{Bgd}}_1+Q\,\beta^{\,\sssize\text{Bgd}}_2,\\
&\gamma^{\sssize\text{Bgd}}_{21}=\partial_y\beta^{\,\sssize\text{Bgd}}_2
-S\,\beta^{\,\sssize\text{Bgd}}_1+R\,\beta^{\,\sssize\text{Bgd}}_2.
\endaligned
\mytag{3.3}
$$
The fourth order expressions are given by the formulas \thetag{2.4} 
in \mycite{5}:  
$$
\hskip -2em
\aligned
&\delta^{\,\sssize\text{Bgd}}_{10}=\partial_x\gamma^{\sssize\text{Bgd}}_{10}
-2\,Q\,\gamma^{\sssize\text{Bgd}}_{10}+P\,(\gamma^{\sssize\text{Bgd}}_{20}
+\gamma^{\sssize\text{Bgd}}_{11})-5\,\alpha^{\sssize\text{Bgd}}_0\,
\beta^{\,\sssize\text{Bgd}}_1,\\
&\delta^{\,\sssize\text{Bgd}}_{20}=\partial_x\gamma^{\sssize\text{Bgd}}_{20}
-R\,\gamma^{\sssize\text{Bgd}}_{10}+P\,\gamma^{\sssize\text{Bgd}}_{21}
-4\,\alpha^{\sssize\text{Bgd}}_1\,\beta^{\,\sssize\text{Bgd}}_1
-\alpha^{\sssize\text{Bgd}}_0\,\beta^{\,\sssize\text{Bgd}}_2,\\
&\delta^{\,\sssize\text{Bgd}}_{30}=\partial_y\gamma^{\sssize\text{Bgd}}_{20}
-S\,\gamma^{\sssize\text{Bgd}}_{10}+Q\,\gamma^{\sssize\text{Bgd}}_{21}
-4\,\alpha^{\sssize\text{Bgd}}_2\,\beta^{\,\sssize\text{Bgd}}_1
-\alpha^{\sssize\text{Bgd}}_1\,\beta^{\,\sssize\text{Bgd}}_2,\\
&\delta^{\,\sssize\text{Bgd}}_{11}=\partial_x\gamma^{\sssize\text{Bgd}}_{11}
-R\,\gamma^{\sssize\text{Bgd}}_{10}+P\,\gamma^{\sssize\text{Bgd}}_{21}
-\alpha^{\sssize\text{Bgd}}_1\,\beta^{\,\sssize\text{Bgd}}_1
-4\,\alpha^{\sssize\text{Bgd}}_0\,\beta^{\,\sssize\text{Bgd}}_2,\\
&\delta^{\,\sssize\text{Bgd}}_{21}=\partial_x\gamma^{\sssize\text{Bgd}}_{21}
-R\,(\gamma^{\sssize\text{Bgd}}_{20}+\gamma^{\sssize\text{Bgd}}_{11})
+2\,Q\,\gamma^{\sssize\text{Bgd}}_{21}-5\,\alpha^{\sssize\text{Bgd}}_1
\,\beta^{\,\sssize\text{Bgd}}_2,\\
&\delta^{\,\sssize\text{Bgd}}_{31}=\partial_y\gamma^{\sssize\text{Bgd}}_{21}
-S\,(\gamma^{\sssize\text{Bgd}}_{20}+\gamma^{\sssize\text{Bgd}}_{11})
+2\,R\,\gamma^{\sssize\text{Bgd}}_{21}-5\,\alpha^{\sssize\text{Bgd}}_2
\,\beta^{\,\sssize\text{Bgd}}_2.
\endaligned
\mytag{3.4}
$$
The fifth order expressions are given by the formulas \thetag{2.5} 
in \mycite{5}:  
$$
\pagebreak
\hskip -2em
\aligned
&\epsilon^{\,\sssize\text{Bgd}}_{10}=\partial_x\delta^{\,\sssize\text{Bgd}}_{10}
-3\,Q\,\delta^{\,\sssize\text{Bgd}}_{10}+P\,(2\,\delta^{\,\sssize\text{Bgd}}_{20}
+\delta^{\,\sssize\text{Bgd}}_{11})-12\,\alpha^{\sssize\text{Bgd}}_0\,
\gamma^{\sssize\text{Bgd}}_{10},\\
&\epsilon^{\,\sssize\text{Bgd}}_{20}=\partial_y\delta^{\,\sssize\text{Bgd}}_{10}
-3\,R\,\delta^{\,\sssize\text{Bgd}}_{10}+Q\,(2\,\delta^{\,\sssize\text{Bgd}}_{20}
+\delta^{\,\sssize\text{Bgd}}_{11})-12\,\alpha^{\sssize\text{Bgd}}_1\,
\gamma^{\sssize\text{Bgd}}_{10},\\
&\aligned
\epsilon^{\,\sssize\text{Bgd}}_{11}&=\partial_x\delta^{\,\sssize\text{Bgd}}_{11}
-R\,\delta^{\,\sssize\text{Bgd}}_{10}-Q\,\delta^{\,\sssize\text{Bgd}}_{11}
+2\,P\,\delta^{\,\sssize\text{Bgd}}_{21}-2\,\alpha^{\sssize\text{Bgd}}_1
\,\gamma^{\sssize\text{Bgd}}_{10}\,-\\
\vspace{-1ex}
&-\,10\,\alpha^{\sssize\text{Bgd}}_0\,\gamma^{\sssize\text{Bgd}}_{11}
-10\,(\beta^{\,\sssize\text{Bgd}}_1)^2.
\endaligned
\endaligned
\mytag{3.5}
$$
And finally, the sixth order expression is given by the formula \thetag{2.6} 
in \mycite{5}:  
$$
\hskip -2em
\lambda^{\,\sssize\text{Bgd}}_{10}=\partial_x\epsilon^{\,\sssize\text{Bgd}}_{10}
-4\,Q\,\epsilon^{\,\sssize\text{Bgd}}_{10}+P\,(3\,\epsilon^{\,\sssize\text{Bgd}}_{20}
+\epsilon^{\,\sssize\text{Bgd}}_{11})-21\,\alpha^{\sssize\text{Bgd}}_0
\,\delta^{\,\sssize\text{Bgd}}_{10}.
\mytag{3.6}
$$\par
     Comparing \mythetag{3.1} with the formulas \thetag{2.15} in \mycite{7}, 
one can easily formulate and prove the following lemma.\par
\mylemma{3.1} Bagderina's alpha quantities \mythetag{3.1} coincide with the 
components of the symmetric two-dimensional array $\Omega$ constructed
in \mycite{7}:
$$
\xalignat 3
&\hskip -2em
\alpha^{\sssize\text{Bgd}}_0=\Omega_{11},
&&\alpha^{\sssize\text{Bgd}}_1=\Omega_{12}=\Omega_{21},
&&\alpha^{\sssize\text{Bgd}}_2=\Omega_{22}. 
\quad
\mytag{3.7}
\endxalignat
$$
\endproclaim
     From \mycite{7} it is known, that the quantities $\Omega_{ij}$ 
in \mythetag{3.7} constitute nether a tensor invariant nor a pseudotensor 
invariant. The transformation rule for them is similar to \mythetag{2.5},
but a little bit different. It is given by the formula \thetag{2.13} in
\mycite{7}.\par
\mylemma{3.2} Bagderina's beta quantities \mythetag{3.2} coincide with the 
components of the pseudocovectorial field\/ $\boldsymbol\alpha$ of the weight 
$1$ constructed in \mycite{7}:
$$
\xalignat 2
&\hskip -2em
\beta^{\,\sssize\text{Bgd}}_1=\alpha_1=A,
&&\beta^{\,\sssize\text{Bgd}}_2=\alpha_2=B.
\mytag{3.8}
\endxalignat
$$
\endproclaim
     Lemma~\mythelemma{3.2} is proved by comparing the formulas \thetag{1.6}
from \mycite{5} with \thetag{2.19} in \mycite{7}.\par
     Apart from $\boldsymbol\alpha$, another pseudotensorial field $\bold d$
was introduced in \mycite{7}. Its components are given by the following 
skew-symmetric matrix in any local coordinates:
$$
\hskip -2em
d_{ij}=\Vmatrix\format \r&\quad\l\\ 0 & 1\\-1 & 0\endVmatrix.
\mytag{3.9}
$$
The type of the field $\bold d$ with the components \mythetag{3.9} is $(0,2)$,
its weight is $-1$. The same matrix \mythetag{3.9} provide the components of 
a pseudotensor of the type $(2,0)$:
$$
\hskip -2em
d^{\kern 1pt ij}=\Vmatrix\format \r&\quad\l\\ 0 & 1\\-1 & 0\endVmatrix.
\mytag{3.10}
$$
This field is denoted by the same symbol $\bold d$, its weight is $1$. The 
fields \mythetag{3.9} and \mythetag{3.10} are used for raising in lowering indices
of pseudotensorial invariants of the equation \mythetag{1.1}. They pay the same role
as the metric tensors in metric geometry. In particular, the formulas 
\mythetag{3.8} can be written as 
$$
\xalignat 2
&\hskip -2em
\beta^{\,\sssize\text{Bgd}}_1=A=-\alpha^2,
&&\beta^{\,\sssize\text{Bgd}}_2=B=\alpha^1,
\mytag{3.11}
\endxalignat
$$
where $\alpha^1$ and $\alpha^2$ are the components of a pseudovectorial field
of the weight $2$:
$$
\hskip -2em
\alpha^i=\sum^2_{k=1}d^{\kern 1pt ik}\,\alpha_k.
\mytag{3.12}
$$ 
This field with the components \mythetag{3.12} in \mycite{7} is denoted by 
the same symbol  $\boldsymbol\alpha$ as the previous field with the
components \mythetag{3.8}.\par 
\mytheorem{3.1} The case of the equations of the ninth type of Bagderina's 
classification in \mycite{5} coincide with the case of {\tenbfit maximal degeneration}
in \myciterange{3}{3}{, }{4}.
\endproclaim    
     According to \mycite{3} (see page 7), the simultaneous vanishing condition 
$A=0$ and $B=0$ determines the case of {\bf maximal degeneration}. In order to
prove Theorem~\mythetheorem{3.1} now it is sufficient to compare this condition 
with the condition $\beta^{\,\sssize\text{Bgd}}_1=0$ and 
$\beta^{\,\sssize\text{Bgd}}_2=0$ in item 9 of Theorem 2 in \mycite{5} and apply 
the above formulas \mythetag{3.8}.\par 
     Apart from \mythetag{3.1}, \mythetag{3.2}, \mythetag{3.3}, \mythetag{3.4},
\mythetag{3.5}, and \mythetag{3.6}, Yu\.~Yu\. Bagderina uses another series
of huge special notations in \mycite{5}. The formulas \thetag{2.16} from
\mycite{5} yield 
$$
\hskip -2em
\aligned
&J^{\,\sssize\text{Bgd}}_0
=(\beta^{\,\sssize\text{Bgd}}_2)^2\,\gamma^{\sssize\text{Bgd}}_{10}
-\beta^{\,\sssize\text{Bgd}}_1\,\beta^{\,\sssize\text{Bgd}}_2\,
(\gamma^{\sssize\text{Bgd}}_{20}+\gamma^{\sssize\text{Bgd}}_{11})
+(\beta^{\,\sssize\text{Bgd}}_1)^2\,\gamma^{\sssize\text{Bgd}}_{21},\\
\vspace{1ex}
&\aligned
J^{\,\sssize\text{Bgd}}_1=
\beta^{\,\sssize\text{Bgd}}_2\,(&\delta^{\,\sssize\text{Bgd}}_{20}
-\delta^{\,\sssize\text{Bgd}}_{11})
+\beta^{\,\sssize\text{Bgd}}_1\,(\delta^{\,\sssize\text{Bgd}}_{21}
-\delta^{\,\sssize\text{Bgd}}_{30})\,+\\
\vspace{-1ex}
&+\,\dfrac{7}{5}\,(\gamma^{\sssize\text{Bgd}}_{20}
-\gamma^{\sssize\text{Bgd}}_{11})^2
+\dfrac{3}{5}\,(\gamma^{\sssize\text{Bgd}}_{20}\,
\gamma^{\sssize\text{Bgd}}_{11}-\gamma^{\sssize\text{Bgd}}_{10}\,
\gamma^{\sssize\text{Bgd}}_{21}),
\endaligned\\
&\aligned
J^{\,\sssize\text{Bgd}}_2&=\Gamma^{\,\sssize\text{Bgd}}_1
\,(\delta^{\,\sssize\text{Bgd}}_{20}-\delta^{\,\sssize\text{Bgd}}_{11})+\Gamma^{\,\sssize\text{Bgd}}_0\,(\delta^{\,\sssize\text{Bgd}}_{21}
-\delta^{\,\sssize\text{Bgd}}_{30})+3\,(\gamma^{\sssize\text{Bgd}}_{20}\,-\\
\vspace{-1ex}
&-\,\gamma^{\sssize\text{Bgd}}_{11})\,(\gamma^{\sssize\text{Bgd}}_{10}\,\gamma^{\sssize\text{Bgd}}_{21}
-\gamma^{\sssize\text{Bgd}}_{20}\,\gamma^{\sssize\text{Bgd}}_{11})
+\frac{4}{3}\,(\gamma^{\sssize\text{Bgd}}_{20}
-\gamma^{\sssize\text{Bgd}}_{11})^3,
\endaligned\\
&\aligned
J^{\,\sssize\text{Bgd}}_3&=
(\beta^{\,\sssize\text{Bgd}}_2)^3\,\delta^{\,\sssize\text{Bgd}}_{10}
-\beta^{\,\sssize\text{Bgd}}_1\,(\beta^{\,\sssize\text{Bgd}}_2)^2\,
(2\,\delta^{\,\sssize\text{Bgd}}_{20}+\delta^{\,\sssize\text{Bgd}}_{11})
\,+\\
&+\,(\beta^{\,\sssize\text{Bgd}}_1)^2\,\beta^{\,\sssize\text{Bgd}}_2\,
(\delta^{\,\sssize\text{Bgd}}_{30}+2\,\delta^{\,\sssize\text{Bgd}}_{21})
-(\beta^{\,\sssize\text{Bgd}}_1)^3\,\delta^{\,\sssize\text{Bgd}}_{31}\,+\\
&+\,4\,(\gamma^{\sssize\text{Bgd}}_{20}
-\gamma^{\sssize\text{Bgd}}_{11})\,J^{\,\sssize\text{Bgd}}_0,
\endaligned\\
&\aligned
J^{\,\sssize\text{Bgd}}_4&=
-\beta^{\,\sssize\text{Bgd}}_2\,(\beta^{\,\sssize\text{Bgd}}_2
\,\Gamma^{\,\sssize\text{Bgd}}_0+2\,\beta^{\,\sssize\text{Bgd}}_1
\,\Gamma^{\,\sssize\text{Bgd}}_1)\,(2\,\delta^{\,\sssize\text{Bgd}}_{20}+\delta^{\,\sssize\text{Bgd}}_{11})\,+\\
&+\,\beta^{\,\sssize\text{Bgd}}_1\,(2\,\beta^{\,\sssize\text{Bgd}}_2
\,\Gamma^{\,\sssize\text{Bgd}}_0
+\beta^{\,\sssize\text{Bgd}}_1\,\Gamma^{\,\sssize\text{Bgd}}_1)
\,(\delta^{\,\sssize\text{Bgd}}_{30}+2\,\delta^{\,\sssize\text{Bgd}}_{21})\,+\\
\vspace{-1ex}
&+3\,(\beta^{\,\sssize\text{Bgd}}_2)^2\,\Gamma^{\,\sssize\text{Bgd}}_1
\,\delta^{\,\sssize\text{Bgd}}_{10}-3\,(\beta^{\,\sssize\text{Bgd}}_1)^2
\,\Gamma^{\,\sssize\text{Bgd}}_0\,\delta^{\,\sssize\text{Bgd}}_{31}
+\frac{66}{5}\,\times\\
\vspace{-1ex}
&\times\,(\gamma^{\sssize\text{Bgd}}_{20}
-\,\gamma^{\sssize\text{Bgd}}_{11})^2\,J^{\,\sssize\text{Bgd}}_0
+\frac{36}{5}\,(\gamma^{\sssize\text{Bgd}}_{10}\,\gamma^{\sssize\text{Bgd}}_{21}
-\gamma^{\sssize\text{Bgd}}_{20}\,\gamma^{\sssize\text{Bgd}}_{11})
\,J^{\,\sssize\text{Bgd}}_0.
\endaligned
\endaligned
\mytag{3.13}
$$
The quantities $\Gamma^{\,\sssize\text{Bgd}}_0$ and $\Gamma^{\,\sssize\text{Bgd}}_1$
from \mythetag{3.13} are given by the formulas \thetag{2.17} in \mycite{5}:  
$$
\hskip -2em
\aligned
&\Gamma^{\,\sssize\text{Bgd}}_0=3\,\beta^{\,\sssize\text{Bgd}}_2\,
\gamma^{\sssize\text{Bgd}}_{10}+\beta^{\,\sssize\text{Bgd}}_1\,
(\gamma^{\sssize\text{Bgd}}_{20}-4\,\gamma^{\sssize\text{Bgd}}_{11}),\\
&\Gamma^{\,\sssize\text{Bgd}}_1=\beta^{\,\sssize\text{Bgd}}_2
\,(4\,\gamma^{\sssize\text{Bgd}}_{20}-\gamma^{\sssize\text{Bgd}}_{11})
-3\,\beta^{\,\sssize\text{Bgd}}_1\,\gamma^{\sssize\text{Bgd}}_{21}.
\endaligned
\mytag{3.14}
$$
\mylemma{3.3} Bagderina's quantity $J^{\,\sssize\text{Bgd}}_0$ from
\mythetag{3.13} is related to the pseudoscalar field $F$ of the weight\/ $1$
constructed in \mycite{7} by means of the formula
$$
\hskip -2em
J^{\,\sssize\text{Bgd}}_0=-F^5.
\mytag{3.15}
$$
\endproclaim 
     Lemma~\mythelemma{3.3} is proved by direct calculations with the use
of the first formula \mythetag{3.13} and the formula \thetag{2.25} from
\mycite{7}.\par
\mytheorem{3.2} The case of the equations of the first type of Bagderina's 
classification in \mycite{5} coincide with the case of {\tenbfit general
position} in \myciterange{3}{3}{, }{4}.
\endproclaim    
     According to \mycite{3} (see section 3 on page 6), the case of {\bf general
position} is defined by the non-vanishing condition for the pseudoscalar invariant
$F$:
$$
\hskip -2em
F\neq 0.
\mytag{3.16}
$$
Due to \mythetag{3.15} in order to prove Theorem~\mythetheorem{3.2} now it is 
sufficient to compare the condition \mythetag{3.16} with the condition 
$J^{\,\sssize\text{Bgd}}_0\neq 0$ in item 1 of Theorem 2 in \mycite{5}.\par 
     The equality \mythetag{3.15} means that Bagderina's quantity 
$J^{\,\sssize\text{Bgd}}_0$ is a pseudoscalar field of the weight $5$. For the 
case of general position Yu\.~Yu\.~Bagderina defines the following quantity in
her paper \mycite{5} (see item 1 of theortem 2 in \mycite{5}):
$$
\hskip -2em
\mu^{\,\sssize\text{Bgd}}_1=(J^{\,\sssize\text{Bgd}}_0)^{1/5}.
\mytag{3.17}
$$
\mylemma{3.4} Bagderina's quantity $\mu^{\,\sssize\text{Bgd}}_1$ from
\mythetag{3.17} is related to the pseudoscalar field $F$ of the weight\/ $1$
constructed in \mycite{7} by means of the formula
$$
\hskip -2em
\mu^{\,\sssize\text{Bgd}}_1=-F.
\mytag{3.18}
$$
\endproclaim 
     The equality \mythetag{3.18} and Lemma~\mythelemma{3.4} are immediate from 
\mythetag{3.15} and \mythetag{3.17}. 
\mylemma{3.5} Bagderina's gamma quantities \mythetag{3.14} coincide with the 
components of the pseudocovectorial field\/ $\boldsymbol\beta$ of the weight 
$3$ constructed in \mycite{7}:
$$
\xalignat 2
&\hskip -2em
\Gamma^{\,\sssize\text{Bgd}}_0=\beta_1=-H,
&&\Gamma^{\,\sssize\text{Bgd}}_1=\beta_2=G.
\mytag{3.19}
\endxalignat
$$
\endproclaim
     The pseudovectorial field\/ $\boldsymbol\beta$ is produced from the 
pseudocovectorial field\/ $\boldsymbol\beta$ in \mycite{7} by raising 
the index of its components with the use of the formula 
$$
\hskip -2em
\beta^{\kern 1pt i}=\sum^2_{k=1}d^{\kern 1pt ik}\,\beta_k.
\mytag{3.20}
$$ 
Then due to \mythetag{3.20} the formulas \mythetag{3.19} are written as 
follows: 
$$
\xalignat 2
&\hskip -2em
\Gamma^{\,\sssize\text{Bgd}}_1=\beta^1=G,
&&\Gamma^{\,\sssize\text{Bgd}}_0=-\beta^2=-H
\mytag{3.21}
\endxalignat
$$
The weight of the pseudovectorial field\/ $\boldsymbol\beta$ produced
by means of the above formula \mythetag{3.20} in \mycite{3} is equal to $4$.
\par
     In her paper \mycite{5} Yu\.~Yu\.~Bagderina defines two differential
operators (see Theo\-rem 2). They are given by the formulas \thetag{2.8} 
in \mycite{5}:
$$
\hskip -2em
\aligned
&\Cal D^{\,\sssize\text{Bgd}}_1=(\mu^{\,\sssize\text{Bgd}}_1)^{-2}
(\beta^{\,\sssize\text{Bgd}}_2\,D_x-\beta^{\,\sssize\text{Bgd}}_1\,D_y),\\
&\Cal D^{\,\sssize\text{Bgd}}_2=\mu^{\,\sssize\text{Bgd}}_2
\,(\beta^{\,\sssize\text{Bgd}}_2\,D_x-\beta^{\,\sssize\text{Bgd}}_1\,D_y)-3\,
(\beta^{\,\sssize\text{Bgd}}_1)^{-1}\,\mu^{\,\sssize\text{Bgd}}_1\,D_x,
\endaligned
\mytag{3.22}
$$
where $D_x$ and $D_y$ can be understood as partial derivatives 
$\partial/\partial x$ and $\partial/\partial y$ acting upon functions 
of the form $f(x,y)$. In the case of {\bf general position} the quantity $\mu^{\,\sssize\text{Bgd}}_2$ is defined by one of the formulas \thetag{2.9} 
in \mycite{5}:
$$
\hskip -2em
\mu^{\,\sssize\text{Bgd}}_2=\Gamma^{\,\sssize\text{Bgd}}_0\,
(\beta^{\,\sssize\text{Bgd}}_1)^{-1}\,(J^{\,\sssize\text{Bgd}}_0)^{-4/5}.
\mytag{3.23}
$$
Applying \mythetag{3.11} and \mythetag{3.18} to the first formula
\mythetag{3.22}, we get
$$
\pagebreak
\hskip -2em
\Cal D^{\,\sssize\text{Bgd}}_1=\frac{\alpha^1}{F^2}
\,\frac{\partial}{\partial x}
+\frac{\alpha^2}{F^2}\,\frac{\partial}{\partial y}.
\mytag{3.24}
$$\par
     The operator $\Cal D^{\,\sssize\text{Bgd}}_2$ is a little bit more
complicated. In order to transform it one should use the following formula
from \mycite{5} (see Remark 1): 
$$
\hskip -2em
3\,J^{\,\sssize\text{Bgd}}_0=\beta^{\,\sssize\text{Bgd}}_2\,
\Gamma^{\,\sssize\text{Bgd}}_0-\beta^{\,\sssize\text{Bgd}}_1\,
\Gamma^{\,\sssize\text{Bgd}}_1.
\mytag{3.25}
$$
The formula \mythetag{3.25} coincides with the formula \thetag{2.24}
in \mycite{7}. Applying this formula in calculating the operator 
$\Cal D^{\,\sssize\text{Bgd}}_2$, we get the formula
$$
\hskip -2em
\Cal D^{\,\sssize\text{Bgd}}_2=(J^{\,\sssize\text{Bgd}}_0)^{-4/5}
\,(\Gamma^{\,\sssize\text{Bgd}}_1\,D_x-\Gamma^{\,\sssize\text{Bgd}}_0\,D_y).
\mytag{3.26}
$$
The formula \mythetag{3.26} is comprised in \mycite{5} (see Remark 1). 
Applying the formulas \mythetag{3.21} and \mythetag{3.15} to the
above formula \mythetag{3.26}, we get
$$
\hskip -2em
\Cal D^{\,\sssize\text{Bgd}}_2=\frac{\beta^1}{F^4}
\,\frac{\partial}{\partial x}
+\frac{\beta^2}{F^4}\,\frac{\partial}{\partial y}.
\mytag{3.27}
$$
\mytheorem{3.3} Bagderina's differential operators \mythetag{3.22} coincide 
with the vectorial fields $\bold X$ and $\bold Y$ constructed in \mycite{7}:
$$
\xalignat 2
&\hskip -2em
\Cal D^{\,\sssize\text{Bgd}}_1=\bold X,
&&\Cal D^{\,\sssize\text{Bgd}}_2=\bold Y.
\mytag{3.28}
\endxalignat
$$
\endproclaim
     The components of the vector fields $\bold X$ and $\bold Y$ are determined
by the formulas \thetag{3.1} in \mycite{7}. Therefore in order to prove 
the relationships \mythetag{3.28} in Theorem~\mythetheorem{3.3} it is sufficient 
to compare the formulas \thetag{3.1} from \mycite{7} with \mythetag{3.24} and 
\mythetag{3.27}.\par
\head
4. Scalar invariants. 
\endhead
     According to Definition~\mythedefinition{2.2}, scalar and pseudoscalar invariants
are those whose type is $(0,0)$. Scalar invariants differ from pseudoscalar ones by
their weight.  The weight of scalar invariants is zero, while  pseudoscalar invariants
have nonzero weights. In \mycite{7} ten scalar invariants are considered. They are 
denoted through $I_1$, $I_2$, $I_3$, $I_4$, $I_5$, $I_6$, $I_7$, $I_8$, $L$ and $K$. 
Not all of them are independent. We have
$$
\xalignat 4
&\hskip -2em
I_1=-4\,I_6,
&&I_2=\frac{1}{3},
&&I_4=4\,I_6,
&&I_5=-I_8.
\quad
\mytag{4.1}
\endxalignat
$$
The formulas \mythetag{4.1} coincide with the formulas \thetag{3.10} in 
\mycite{4}. Apart from the formulas \mythetag{4.1} we have the following 
relationships:
$$
\xalignat 2
&\hskip -2em
I_5=I_3-L,
&&I_6=-I_1+K.
\mytag{4.2}
\endxalignat
$$
The formulas \mythetag{4.2} are taken from \thetag{3.13} and \thetag{3.14} 
in \mycite{7}. From \mythetag{4.1} and \mythetag{4.2} we derive the following 
expressions for $L$ and $K$:
$$
\xalignat 2
&\hskip -2em
L=I_3+I_8,
&&K=-3\,I_6.
\mytag{4.3}
\endxalignat
$$
Due to \mythetag{4.1} and \mythetag{4.3} all of the ten invariants $I_1$, 
$I_2$, $I_3$, $I_4$, $I_5$, $I_6$, $I_7$, $I_8$, $L$, $K$ are expressed 
through four of them: $I_3$, $I_6$, $I_7$, $I_8$.\par
     The invariants $I_3$, $I_6$, $I_7$, $I_8$ are given by explicit formulas 
taken from \mycite{7} and \mycite{4}. The invariant $I_3$ is given by the following
formula:
$$
\gathered
I_3=\frac{B\,(H\,G_{\sssize 1.0}-G\,H_{\sssize 1.0})}{3\,F^9}
-\frac{A\,(H\,G_{\sssize 0.1}-G\,H_{\sssize 0.1})}{3\,F^9}+
\frac{H\,F_{\sssize 0.1}+G\,F_{\sssize 1.0}}{3\,F^5}\,+\\
\vspace{1ex}
+\,\frac{B\,G^2\,P}{3\,F^9}-\frac{(A\,G^2-2\,H\,B\,G)\,Q}{3\,F^9}+
\frac{(B\,H^2-2\,H\,A\,G)\,R}{3\,F^9}-\frac{A\,H^2\,S}{3\,F^9}.
\endgathered
\quad
\mytag{4.4}
$$
The formula \mythetag{4.4} coincides with the formula \thetag{3.8} in
\mycite{7} and the formula \thetag{3.6} in \mycite{4}. The quantities $A$ and
$B$ in \mythetag{4.4} are given by the formulas \mythetag{3.8} and 
\mythetag{3.2}. The quantities $G$ and $H$ are given by the formulas 
\mythetag{3.21} and \mythetag{3.14}. The quantity $F$ is given by
the formulas \mythetag{3.15} and \mythetag{3.18}. Here are more explicit 
formulas for all of these five quantities $A$, $B$, $G$, $H$, and $F$:
$$
\gather
\hskip -2em
\aligned
&\aligned
 A=P_{\sssize 0.2}&-2\,Q_{\sssize 1.1}+R_{\sssize 2.0}+
 2\,P\,S_{\sssize 1.0}+S\,P_{\sssize 1.0}-\\
 \vspace{0.5ex}
 &-3\,P\,R_{\sssize 0.1}-3\,R\,P_{\sssize 0.1}-
 3\,Q\,R_{\sssize 1.0}+6\,Q\,Q_{\sssize 0.1},
 \endaligned\\
 \vspace{1ex}
&\aligned
 B=S_{\sssize 2.0}&-2\,R_{\sssize 1.1}+Q_{\sssize 0.2}-
 2\,S\,P_{\sssize 0.1}-P\,S_{\sssize 0.1}+\\
 \vspace{0.5ex}
 &+3\,S\,Q_{\sssize 1.0}+3\,Q\,S_{\sssize 1.0}+
 3\,R\,Q_{\sssize 0.1}-6\,R\,R_{\sssize 1.0},
 \endaligned
\endaligned
\quad
\mytag{4.5}\\
\vspace{2ex}
\hskip -2em
\aligned
G&=-B\,B_{\sssize 1.0}-3\,A\,B_{\sssize 0.1}+4\,B\,
A_{\sssize 0.1}+3\,S\,A^2-6\,R\,B\,A+3\,Q\,B^2,\\
\vspace{1ex}
H&=-A\,A_{\sssize 0.1}-3\,B\,A_{\sssize 1.0}+4\,A\,
B_{\sssize 1.0}-3\,P\,B^2+6\,Q\,A\,B-3\,R\,A^2,
\endaligned
\quad
\mytag{4.6}\\
\vspace{2ex}
\hskip -2em
\aligned
F^5=A\,B\,A_{\sssize 0.1}&+B\,A\,B_{\sssize 1.0}-
A^2\,B_{\sssize 0.1}-B^2\,A_{\sssize 1.0}-\\
&-P\,B^3+3\,Q\,A\,B^2-3\,R\,A^2\,B+S\,A^3.
\endaligned
\quad
\mytag{4.7}
\endgather
$$
The formulas \mythetag{4.5} coincide with \thetag{2.19} in \mycite{7}
and with \thetag{1.6} in \mycite{4}. Similarly, the formulas \mythetag{4.6} 
coincide with \thetag{2.23} in \mycite{7} and with \thetag{2.4} in \mycite{4}.
The formula \mythetag{4.7} coincides with \thetag{2.25} in \mycite{7} and
\thetag{2.2} in \mycite{4}.\par 
    Like in \mythetag{2.2}, \mythetag{2.3}, \mythetag{2.4}, and \mythetag{3.1}, 
double indices in \mythetag{4.4}, \mythetag{4.5}, \mythetag{4.6}, and 
\mythetag{4.7} are used to denote partial derivatives. For a given function 
$f(x,y)$ we write
$$
\hskip -2em
f_{\ssize p.q}=\frac{\partial^{p+q}f}{\partial x^p\,\partial y^q}.
\mytag{4.8}
$$
\par
     The invariant $I_6$ is given by the formula \thetag{3.7} in \mycite{4}.
Unfortunately this formula is mistyped. The sign of the last term should be
altered. The correct formula is 
$$
\aligned
I_6&=\frac{A\,(G\,A_{\sssize 0.1}+H\,B_{\sssize 0.1})}{12\,F^7}
-\frac{B\,(G\,A_{\sssize 1.0}+H\,B_{\sssize 1.0})}{12\,F^7}
-\frac{4\,(A\,F_{\sssize 0.1}-B\,F_{\sssize 1.0})}{12\,F^3}\,-\\
\vspace{1ex}
&-\frac{G\,B^2\,P}{12\,F^7}-\frac{(H\,B^2-2\,G\,B\,A)\,Q}{12\,F^7}
-\frac{(G\,A^2-2\,H\,B\,A)\,R}{12\,F^7}-\frac{H\,A^2\,S}{12\,F^7}.
\endaligned
\quad
\mytag{4.9}
$$
Partial derivatives in \mythetag{4.9} are denoted according to the 
convention \mythetag{4.8}. There is another formula for the invariant 
$I_6$ in \mycite{4}, which is more simple:
$$
\hskip -2em
I_6=\frac{A_{\sssize 0.1}-B_{\sssize 1.0}}{3\,F^2}-
\frac{A\,F_{\sssize 0.1}-B\,F_{\sssize 1.0}}{3\,F^3}.
\mytag{4.10}
$$
The formula \mythetag{4.10} coincides with \thetag{3.11} in \mycite{4}. Both
formulas \mythetag{4.9} and \mythetag{4.10} give the same result, though they
look different.\par
     The invariants $I_7$ and $I_8$ are given by the following two formulas:
$$
\gather
\hskip -2em
\aligned
I_7=&\frac{G\,H\,G_{\sssize 1.0}-
G^2\,H_{\sssize 1.0}+H^2\,G_{\sssize 0.1}
-H\,G\,H_{\sssize 0.1}}{3\,F^{11}}\,+\\
\vspace{1ex}
&\hskip 5em+\frac{G^3\,P+3\,G^2\,H\,\,Q+3\,G\,H^2\,R+H^3\,S}
{3\,F^{11}},
\endaligned
\mytag{4.11}\\
\vspace{3ex}
\hskip -2em
\aligned
I_8&=\frac{G\,(A\,G_{\sssize 1.0}+B\,H_{\sssize 1.0})}{3\,F^9}
+\frac{H\,(A\,G_{\sssize 0.1}+B\,H_{\sssize 0.1})}{3\,F^9}
-\frac{10\,(H\,F_{\sssize 0.1}+G\,F_{\sssize 1.0})}{3\,F^5}\,-\\
\vspace{1ex}
&-\frac{B\,G^2\,P}{3\,F^9}+\frac{(A\,G^2-2\,H\,B\,G)\,Q}{3\,F^9}-
\frac{(B\,H^2-2\,H\,A\,G)\,R}{3\,F^9}+\frac{A\,H^2\,S}{3\,F^9}.
\endaligned
\quad
\mytag{4.12}
\endgather
$$
The formula \mythetag{4.11} coincides with the formula \thetag{3.9} in
\mycite{7} and with the formula \thetag{3.8} in \mycite{4}. The formula
\mythetag{4.12} is taken from \thetag{3.9} in \mycite{4}.\par
     Yu\.~Yu\.~Bagderina presented her own invariants in \mycite{5} . 
In the case of {\bf general position} they are given by the formulas
$$
\xalignat 2
&\hskip -2em
I^{\,\sssize\text{Bgd}}_1=\frac{J^{\,\sssize\text{Bgd}}_1}
{(J^{\,\sssize\text{Bgd}}_0)^{4/5}},
&&I^{\,\sssize\text{Bgd}}_2=\frac{J^{\,\sssize\text{Bgd}}_2}
{(J^{\,\sssize\text{Bgd}}_0)^{6/5}},\\
\vspace{-1ex}
\mytag{4.13}\\
\vspace{-1ex}
&\hskip -2em
I^{\,\sssize\text{Bgd}}_3=\frac{J^{\,\sssize\text{Bgd}}_3}
{(J^{\,\sssize\text{Bgd}}_0)^{7/5}},
&&I^{\,\sssize\text{Bgd}}_4=\frac{J^{\,\sssize\text{Bgd}}_4}
{(J^{\,\sssize\text{Bgd}}_0)^{9/5}},\\
\endxalignat
$$
where $J^{\,\sssize\text{Bgd}}_0$, $J^{\,\sssize\text{Bgd}}_1$, 
$J^{\,\sssize\text{Bgd}}_2$, $J^{\,\sssize\text{Bgd}}_3$, 
$J^{\,\sssize\text{Bgd}}_0$ are taken from \mythetag{3.13}
(see Theorem 2 in \mycite{5}). On page 3 of her paper \mycite{5}
Yu\.~Yu\.~Bagderina give the comparison formulas relating her
invariants with $I_3$, $I_6$, $I_7$, and $I_8$. In order to verify
these comparison formulas we need to introduce some special 
coordinates. 
\head
5. Special coordinates. 
\endhead
     For the case of {\bf general position} in \mycite{7} two vectorial
fields $\bold X$ and $\bold Y$ were constructed. Their components are 
given by the formulas
$$
\xalignat 2
&\hskip -2em
X^1=\frac{\alpha^1}{F^2}=\frac{B}{F^2},
&&X^2=\frac{\alpha^2}{F^2}=-\frac{A}{F^2},
\mytag{5.1}\\
\vspace{1ex}
&\hskip -2em
Y^1=\frac{\beta^1}{F^4}=\frac{G}{F^4},
&&Y^2=\frac{\beta^2}{F^4}=\frac{H}{F^4}. 
\mytag{5.2}
\endxalignat
$$
According to Theorem~\mythetheorem{3.3} Bagderina's differential operators
$\Cal D^{\,\sssize\text{Bgd}}_1$ and $\Cal D^{\,\sssize\text{Bgd}}_2$ coincide
with the vector fields $\bold X$ and $\bold Y$ whose components are given in
\mythetag{5.1} and \mythetag{5.2}. The commutator of the vector fields 
$\bold X$ and $\bold Y$ is given by the formula 
$$
\hskip -2em
[\bold X,\,\bold Y]=L\,\bold X-K\,\bold Y,
\mytag{5.3}
$$
where $L$ and $K$ are scalar invariants from \mythetag{4.3}. The formula
\mythetag{5.3} coincides with \thetag{3.15} in \mycite{7}. There it is used
in order to define the invariants $L$ and $K$.\par
\mytheorem{5.1} For any two vector field obeying the relationship 
\mythetag{5.3} there are two scalar functions $u$ and $v$ such that 
the vector fields 
$$
\xalignat 2
&\hskip -2em
\tilde\bold X=\frac{\bold X}{u},
&&\tilde\bold Y=\frac{\bold Y}{v},
\mytag{5.4}
\endxalignat
$$
do commute with each other, i\.\,e\. their commutator is zero:
$[\tilde\bold X,\,\tilde\bold Y]=0$.
\endproclaim
    Theorem~\mythetheorem{5.1} is a rather well-known result. It can be considered
as a two-dimensional reduction of the well-known Frobenius theorem on integrability
of involutive distributions (see Proposition 1.2 in Chapter I of \mycite{10}). The
proof of Theorem~\mythetheorem{5.1} is sketched on page 8 of \mycite{7}.\par
    Two commuting vector fields \mythetag{5.4} on the plane define a local 
curvilinear coordinate system $(\tilde x,\tilde y)$. The coordinates $\tilde x$ and
$\tilde y$ are related to the initial coordinates $x$ and $y$ by means of some
transformation \mythetag{1.2}. They are special coordinates for the equation
\mythetag{1.1} since they are related to the vectorial invariants of this equation.
Transforming to these special coordinates $\tilde x$ and $\tilde y$, we would straighten 
the vector fields \mythetag{5.4}, i\.\,e\. we would have the following formulas for 
their components:
$$
\xalignat 2
&\hskip -2em
\tilde X^1=1,
&&\tilde X^2=0,
\mytag{5.5}\\
\vspace{1ex}
&\hskip -2em
\tilde Y^1=0,
&&\tilde Y^2=1. 
\mytag{5.6}
\endxalignat
$$
Applying \mythetag{5.4} to \mythetag{5.5} and \mythetag{5.6}, we derive
$$
\xalignat 2
&\hskip -2em
X^1=u,
&&X^2=0,
\mytag{5.7}\\
\vspace{1ex}
&\hskip -2em
Y^1=0,
&&Y^2=v. 
\mytag{5.8}
\endxalignat
$$
\par
     Now we apply \mythetag{5.7} and \mythetag{5.8} to \mythetag{5.3}. As a 
result we derive 
$$
\xalignat 2
&\hskip -2em
\bold Xv=-K\,v,
&&\bold Yu=-L\,u.
\mytag{5.9}
\endxalignat
$$     
The formulas \mythetag{5.9} coincide with \thetag{3.17} in \mycite{7}. Being
combined with \mythetag{5.7} and \mythetag{5.8}, the formulas \mythetag{5.9}
lead to the differential equations
$$
\xalignat 2
&\hskip -2em
u\,v_{\sssize 1.0}=-K\,v,
&&v\,u_{\sssize 0.1}=-L\,u.
\mytag{5.10}
\endxalignat
$$
The partial derivatives in \mythetag{5.10} refer to the special coordinates
$\tilde x$ and $\tilde y$. However below we shall omit the tilde sign, e\.\,g\.
assuming that these special coordinates were chosen from the very beginning. 
\par
     The next step is to apply \mythetag{5.1} and \mythetag{5.2} to 
\mythetag{5.7} and \mythetag{5.8}. This yields 
$$
\xalignat 4
&\hskip -2em
u=\frac{B}{F^2},
&&A=0,
&&G=0,
&&v=\frac{H}{F^4}.
\quad
\mytag{5.11}
\endxalignat
$$
Taking $A=0$ from \mythetag{5.11} and applying it to \mythetag{4.6} and 
\mythetag{4.7}, we get
$$
\xalignat 3
&\hskip -2em
P=-\frac{F^5}{B^3},
&&H=-3\,P\,B^2,
&&Q=\frac{B_{\sssize 1.0}}{3\,B}.
\quad
\mytag{5.12}
\endxalignat
$$
Combining the formulas \mythetag{5.11} and \mythetag{5.12}, we can write
$$
\xalignat 4
&\hskip -2em
u=\frac{B}{F^2},
&&v=\frac{3\,F}{B},
&&P=-\frac{F^5}{B^3},
&&Q=\frac{B_{\sssize 1.0}}{3\,B}.
\quad
\mytag{5.13}
\endxalignat
$$
The formulas \mythetag{5.13} do coincide with \thetag{3.20} in \mycite{7}. We
use the first two of them in order to substitute them into \mythetag{5.10}. As 
a result we get  
$$
\xalignat 2
&\hskip -2em
L=\frac{6\,B\,F_{\sssize 0.1}-3\,B_{\sssize 0.1}\,F}{B^2},
&&K=\frac{F\,B_{\sssize 1.0}-B\,F_{\sssize 1.0}}{F^3}.
\mytag{5.14}
\endxalignat
$$
The formulas \mythetag{5.14} coincide with the formulas \thetag{3.22} in 
\mycite{7}.\par
     Note that $L$ and $K$ in \mythetag{5.14} are scalar invariants considered 
above in section 4. The other invariants considered in section 4 are $I_1$, $I_2$, 
$I_3$, $I_4$, $I_5$, $I_6$, $I_7$, $I_8$. In \mycite{7} they are introduced through 
the following formulas: 
$$
\xalignat 2
&\hskip -2em
\nabla_{\bold X}\bold X=I_1\,\bold X+I_2\,\bold Y,
&&\nabla_{\bold X}\bold Y=I_3\,\bold X+I_4\,\bold Y,\\
\vspace{-1.7ex}
&&&\mytag{5.15}\\
\vspace{-1.7ex}
&\hskip -2em
\nabla_{\bold Y}\bold X=I_5\,\bold X+I_6\,\bold Y,
&&\nabla_{\bold Y}\bold Y=I_7\,\bold X+I_8\,\bold Y
\endxalignat
$$
(see the formulas \thetag{3.5} in \mycite{7}). The covariant derivatives in the 
left hand sides of the formulas \mythetag{5.15} are calculated with the use of 
the connection components $\varGamma^k_{ij}$ introduced by the formulas 
\thetag{3.4} in \mycite{7}. Irrespective to the choice of a coordinate system
these connection components are given by the formulas
$$
\xalignat 2
&\hskip -2em
\varGamma^1_{11}=Q+\frac{2\,F_{\sssize 1.0}}{3\,F},
&&\varGamma^2_{11}=-P,\\
\vspace{1ex}
&\hskip -2em
\varGamma^1_{12}=\varGamma^1_{21}=R+\frac{F_{\sssize 0.1}}{3\,F},
&&\varGamma^2_{12}=\varGamma^2_{21}=-Q+\frac{F_{\sssize 1.0}}{3\,F},
\mytag{5.16}\\
\vspace{1ex}
&\hskip -2em
\varGamma^1_{22}=S,
&&\varGamma^2_{22}=-R+\frac{2\,F_{\sssize 0.1}}{3\,F}.
\endxalignat
$$
Applying \mythetag{5.7}, \mythetag{5.8} and \mythetag{5.16} to \mythetag{5.15},
one can easily calculate the invariants $I_1$, $I_2$, $I_3$, $I_4$, $I_5$, 
$I_6$, $I_7$, $I_8$ in special coordinates explicitly:
$$
\xalignat 2
&\hskip -2em
I_2=\frac{1}{3},
&&I_1=\frac{4\,F\,B_{\sssize 1.0}-4\,B\,F_{\sssize 1.0}}{3\,F^3},\\
\vspace{1ex}
&\hskip -2em
I_3=\frac{F_{\sssize 0.1}+3\,F\,B}{B},
&&I_4=\frac{4\,B\,F_{\sssize 1.0}-4\,F\,B_{\sssize 1.0}}{3\,F^3},\\
\vspace{-1ex}
\mytag{5.17}\\
\vspace{-1ex}
&\hskip -2em
I_6=\frac{B\,F_{\sssize 1.0}-F\,B_{\sssize 1.0}}{3\,F^3},
&&I_5=\frac{3\,F\,B_{\sssize 0.1}+3\,F\,B\,R-5\,B\,F_{\sssize 0.1}}{B^2},\\
\vspace{1ex}
&\hskip -2em
I_7=\frac{9\,F^4\,S}{B^3},
&&I_8=\frac{5\,B\,F_{\sssize 0.1}-3\,F\,B_{\sssize 0.1}-3\,F\,B\,R}{B^2}.
\endxalignat
$$
The formulas \mythetag{4.1}, \mythetag{4.2}, \mythetag{4.3} are easily derived
from \mythetag{5.14} and \mythetag{5.17}. It is that very way they were derived
in \mycite{7}.\par  
     Note that the formulas \mythetag{5.13} express $P$ and $Q$ through 
$F$ and $B$. Differentiating these expressions, we can calculate partial derivatives 
of $P$ and $Q$, \pagebreak expressing them through partial derivatives of $F$ 
and $B$. Here are the formulas for $P_{\sssize 1.0}$ and $P_{\sssize 0.1}$:
$$
\xalignat 2
&\hskip -2em
P_{\sssize 1.0}=\frac{3\,F^5\,B_{\sssize 1.0}-5\,B\,F^4\,F_{\sssize 1.0}}{B^4},
&&P_{\sssize 0.1}=\frac{3\,F^5\,B_{\sssize 0.1}-5\,B\,F^4\,F_{\sssize 0.1}}{B^4}.
\quad
\mytag{5.18}
\endxalignat
$$
For the partial derivatives $Q_{\sssize 1.0}$ and $Q_{\sssize 0.1}$ we have 
$$
\xalignat 2
&\hskip -2em
Q_{\sssize 1.0}=\frac{B\,B_{\sssize 2.0}-B_{\sssize 1.0}^2}{3\,B^2},
&&Q_{\sssize 0.1}=\frac{B\,B_{\sssize 1.1}
-B_{\sssize 1.0}\,B_{\sssize 0.1}}{3\,B^2}.
\quad
\mytag{5.19}
\endxalignat
$$\par
     Differentiating \mythetag{5.18} and \mythetag{5.19}, we can derive 
formulas for all partial derivatives of $P$ and $Q$ in terms of $F$and $B$ 
and in terms of partial derivatives of them:
$$
\gather
P_{\ssize m.n}=P_{\ssize m.n}(F,B,\ldots,F_{\ssize p.q},B_{\ssize p.q},\,\ldots)
\text{, \ where \ }p+q\leqslant m+n,
\mytag{5.20}\\
\vspace{1ex}
Q_{\ssize m.n}=Q_{\ssize m.n}(B,\ldots,B_{\ssize p.q},\ldots)
\text{, \ where \ }p+q\leqslant m+n.
\mytag{5.21}
\endgather
$$
In particular we have the following formulas:
$$
\align
&\hskip -2em
\aligned
P_{\sssize 0.2}=\frac{3\,F^5\,B_{\sssize 0.2}}{B^4}
&-\frac{5\,F^4\,F_{\sssize 0.2}}{B^3}
-\frac{12\,F^5\,B_{\sssize 0.1}^2}{B^5}\,+\\
&\quad+\,\frac{30\,F^4\,F_{\sssize 0.1}\,B_{\sssize 0.1}}{B^4}
-\frac{20\,F^3\,F_{\sssize 0.1}^2}{B^3},
\endaligned
\mytag{5.22}\\
\vspace{1.5ex}
&\hskip -2em
Q_{\sssize 1.1}=\frac{B_{\sssize 2.1}}{3\,B}
-\frac{B_{\sssize 0.1}\,B_{\sssize 2.0}}{3\,B}
-\frac{2\,B_{\sssize 1.0}\,B_{\sssize 1.1}}{3\,B}
+\frac{2\,B_{\sssize 0.1}\,B_{\sssize 1.0}^2}{3\,B},
\mytag{5.23}\\
\vspace{1.5ex}
&\hskip -2em
Q_{\sssize 0.2}=\frac{B_{\sssize 1.2}}{3\,B}
-\frac{B_{\sssize 1.0}\,B_{\sssize 0.2}}{3\,B}
-\frac{2\,B_{\sssize 0.1}\,B_{\sssize 1.1}}{3\,B}
+\frac{2\,B_{\sssize 1.0}\,B_{\sssize 0.1}^2}{3\,B},
\mytag{5.24}
\endalign
$$
Substituting \mythetag{5.22}, \mythetag{5.23}, \mythetag{5.24} with
\mythetag{5.18} and \mythetag{5.19} into \mythetag{4.5} and taking 
into account $A=0$ from \mythetag{5.11}, we derive the following 
expressions for $R_{\sssize 2.0}$ and $R_{\sssize 1.1}$:
$$
\gather
\hskip -2em
\gathered
R_{\sssize 2.0}=\frac{B_{\sssize 1.0}}{B}\,R_{\sssize 1.0}
-\frac{3\,F^5}{B^3}\,R_{\sssize 0.1}+\biggl(\frac{9\,F^5\,B_{\sssize 0.1}}{B^4}
-\frac{15\,F^4\,F_{\sssize 0.1}}{B^3}\biggr)R\,+\\
\vspace{1ex}
+\frac{2\,F^5}{B^3}\,S_{\sssize 1.0}+\biggl(\frac{5\,F^4\,F_{\sssize 1.0}}{B^3}
-\frac{3\,F^5\,B_{\sssize 1.0}}{B^4}\biggr)S+\frac{2\,B_{\sssize 2.1}}{3\,B}
-\frac{2\,B_{\sssize 0.1}\,B_{\sssize 2.0}}{3\,B^2}\,-\\
\vspace{1ex}
-\frac{2\,B_{\sssize 1.0}\,B_{\sssize 1.1}}{B^2}
-\frac{3\,F^5\,B_{\sssize 0.2}}{B^4}+\frac{5\,F^4\,F_{\sssize 0.2}}{B^3}
+\frac{2\,B_{\sssize 0.1}\,B_{\sssize 1.0}^2}{B^3}
+\frac{20\,F^3\,F_{\sssize 0.1}^2}{B^3}\,-\\
\vspace{1ex}
-\frac{30\,F^4\,F_{\sssize 0.1}\,B_{\sssize 0.1}}{B^4}
+\frac{12\,F^5\,B_{\sssize 0.1}^2}{B^5},
\endgathered
\quad
\mytag{5.25}\\
\vspace{1ex}
\hskip -2em
\gathered
R_{\sssize 1.1}=\frac{1}{2}\,S_{\sssize 2.0}-3\,R\,R_{\sssize 1.0}
+\frac{B_{\sssize 1.0}}{2\,B}\,S_{\sssize 1.0}
+\frac{F^5}{2\,B^3}\,S_{\sssize 0.1}+\biggl(\frac{B_{\sssize 1.0}}{2\,B}\,-
\\
\vspace{1ex}
-\,\frac{B_{\sssize 0.1}\,B_{\sssize 1.0}}{2\,B^2}\biggr)R
+\biggl(\frac{B_{\sssize 2.0}}{2\,B}-\frac{B_{\sssize 1.0}^2}{2\,B^2}
-\frac{3\,F^5\,B_{\sssize 0.1}}{B^4}
+\frac{5\,F^4\,F_{\sssize 0.1}}{B^3}\biggr)S\,+\\
\vspace{1ex}
+\,\frac{B_{\sssize 1.2}}{6\,B}
-\frac{B_{\sssize 0.1}\,B_{\sssize 1.1}}{3\,B^2}
-\frac{B_{\sssize 1.0}\,B_{\sssize 0.2}}{6\,B^2}
+\frac{B_{\sssize 1.0}\,B_{\sssize 0.1}^2}{3\,B^2}
-\frac{B}{2}.
\endgathered
\mytag{5.26}
\endgather
$$
\par
     Differentiating the relationships \mythetag{5.25} and \mythetag{5.26}
we produce the expressions for higher order derivatives of $R$. Some of such
expressions can be produced in two different ways, e\.\,g\. for the derivative 
$R_{\sssize 2.1}$ we have
$$
\xalignat 2
&\hskip -2em
R_{\sssize 2.1}=\frac{\partial R_{\sssize 2.0}}{\partial y},
&&R_{\sssize 2.1}=\frac{\partial R_{\sssize 1.1}}{\partial x}.
\mytag{5.27}
\endxalignat
$$
The right hand sides of \mythetag{5.27} derived from \mythetag{5.25} and 
\mythetag{5.26} are not identically equal to each other. As a result we 
get a non-trivial equation:

$$
\hskip -2em
\frac{\partial R_{\sssize 2.0}}{\partial y}
-\frac{\partial R_{\sssize 1.1}}{\partial x}=0.
\mytag{5.28}
$$
Equations like \mythetag{5.28} are called compatibility conditions in the
theory of differential equations. In our particular case it turns out that
the equation \mythetag{5.28} can be resolved with respect to the derivative
$S_{\sssize 3.0}$.\par
     Higher order partial derivatives of $P$ and $Q$ are calculated by
differentiating \mythetag{5.18} and \mythetag{5.19} (see \mythetag{5.20}
and \mythetag{5.21}). They produce no compatibility conditions like 
\mythetag{5.28}, i\.\,e\. their compatibility conditions are identically 
fulfilled.\par
     The fourth order partial derivatives $S_{\sssize 4.0}$ and $S_{\sssize 3.1}$
are calculated by differentiating the expression for $S_{\sssize 3.0}$ which is
derived from \mythetag{5.28}.\par
     The third order partial derivatives $R_{\sssize 3.0}$, $R_{\sssize 2.1}$,
$R_{\sssize 1.2}$ and the fourth order partial derivatives $R_{\sssize 4.0}$, 
$R_{\sssize 3.1}$, $R_{\sssize 2.2}$, $R_{\sssize 1.3}$ are calculated by
differentiating \mythetag{5.25} and \mythetag{5.26}. They produce no 
compatibility conditions of the order 3 and 4 other than \mythetag{5.28}.\par
     We do not present the results of calculations mentioned just above 
since they are very huge formulas. Upon completing these calculations one
can proceed to Bagderina's formulas \mythetag{3.1}, \mythetag{3.2}, 
\mythetag{3.3}, \mythetag{3.4}. In the special coordinates introduced in
\mycite{7} and used in the present section the formulas \mythetag{3.1} yield
$$
\hskip -2em
\gathered
\alpha^{\sssize\text{Bgd}}_0=\frac{B_{\sssize 2.0}}{3\,B}
-\frac{5\,B_{\sssize 1.0}^2}{9\,B^2}
-\frac{3\,F^5\,B_{\sssize 0.1}}{B^4}
+\frac{5\,F^4\,F_{\sssize 0.1}}{B^3}
-\frac{2\,F^5}{B^3}\,R,\\
\vspace{1ex}
\alpha^{\sssize\text{Bgd}}_1=-\frac{B_{\sssize 1.1}}{3\,B}
+\frac{B_{\sssize 0.1}\,B_{\sssize 1.0}}{3\,B^2}
+R_{\sssize 1.0}
-\frac{B_{\sssize 1.0}}{3\,B}\,R
-\frac{F^5}{B^3}\,S,\\
\vspace{1ex}
\alpha^{\sssize\text{Bgd}}_2=S_{\sssize 1.0}-R_{\sssize 0.1}-2\,R^2
+\frac{2\,B_{\sssize 1.0}}{3\,B}\,S.
\endgathered
\mytag{5.29}
$$
The formulas \mythetag{3.2} in our special coordinates simplify to
$$
\xalignat 2
&\hskip -2em
\beta^{\,\sssize\text{Bgd}}_1=0,
&&\beta^{\,\sssize\text{Bgd}}_2=B,
\mytag{5.30}
\endxalignat
$$   
which is not surprising due to Lemma~\mythelemma{3.2} (see \mythetag{3.8}
and \mythetag{5.11}). The formulas \mythetag{3.3} in our special coordinates 
also yield very simple results:
$$
\xalignat 2
&\hskip -2em
\gamma^{\sssize\text{Bgd}}_{10}=-\frac{F^5}{B^2},
&&\gamma^{\sssize\text{Bgd}}_{11}=\frac{4}{3}\,B_{\sssize 1.0},\\
\vspace{-1.2ex}
\mytag{5.31}\\
\vspace{-1.2ex}
&\hskip -2em
\gamma^{\sssize\text{Bgd}}_{20}=\frac{1}{3}\,B_{\sssize 1.0},
&&\gamma^{\sssize\text{Bgd}}_{21}=B_{\sssize 0.1}+B\,R.
\endxalignat
$$
The results of the formulas \mythetag{3.4} \pagebreak in our special coordinates 
are not very simple. However, they are not very complicated too. Here are these
results: 
$$
\gathered
\delta^{\,\sssize\text{Bgd}}_{10}=\frac{F^5\,B_{\sssize 1.0}}{B^3}
-\frac{5\,F^4\,F_{\sssize 1.0}}{B^2},\\
\vspace{1ex}
\delta^{\,\sssize\text{Bgd}}_{20}=\frac{5\,B_{\sssize 1.0}^2}{9\,B}
+\frac{2\,F^5\,B_{\sssize 0.1}}{B^3}
-\frac{5\,F^4\,F_{\sssize 0.1}}{B^2}
+\frac{2\,F^5}{B^2}\,R,\\
\vspace{1ex}
\delta^{\,\sssize\text{Bgd}}_{30}=\frac{2}{3}\,B_{\sssize 1.1}
-B\,R_{\sssize 1.0}+\frac{2\,B_{\sssize 1.0}}{3}\,R
+\frac{2\,F^5}{B^2}\,S,\\
\vspace{1ex}
\delta^{\,\sssize\text{Bgd}}_{11}=\frac{20\,B_{\sssize 1.0}^2}{9\,B}
+\frac{11\,F^5\,B_{\sssize 0.1}}{B^3}
-\frac{20\,F^4\,F_{\sssize 0.1}}{B^2}
+\frac{8\,F^5}{B^2}\,R,\\
\vspace{1ex}
\delta^{\,\sssize\text{Bgd}}_{21}=\frac{8}{3}\,B_{\sssize 1.1}
-\frac{B_{\sssize 0.1}\,B_{\sssize 1.0}}{B}
-4\,B\,R_{\sssize 1.0}
+\frac{5\,B_{\sssize 1.0}}{3}\,R
+\frac{5\,F^5}{B^2}\,S,\\
\vspace{1ex}
\delta^{\,\sssize\text{Bgd}}_{31}=B_{\sssize 0.2}
+6\,B\,R_{\sssize 0.1}
-5\,B\,S_{\sssize 1.0}
+12\,B\,R^2
+3\,B_{\sssize 0.1}\,R
-5\,B_{\sssize 1.0}\,S.
\endgathered
\mytag{5.32}
$$
The formulas \mythetag{3.5} and \mythetag{3.6} are not used in item 1
of Bagderina's Theorem 2 in \mycite{5}, i\.\,e\. they do not refer to
the case of {\bf general position}. For this reason we shall not 
calculate their values here and, having completed \mythetag{5.29},
\mythetag{5.30}, \mythetag{5.31}, \mythetag{5.32}, we proceed to
the quantities \mythetag{3.14}. The values of these quantities in our
special coordinates are given by the following very simple formulas:
$$
\xalignat 2
&\hskip -2em
\Gamma^{\,\sssize\text{Bgd}}_0=-\frac{3\,F^5}{B},
&&\Gamma^{\,\sssize\text{Bgd}}_1=0,
\mytag{5.33}
\endxalignat
$$   
The formulas \mythetag{5.33} agree with Lemma~\mythelemma{3.5} 
(see \mythetag{3.19}, \mythetag{5.11}, \mythetag{5.12}).\par
     Now let's proceed to the quantities \mythetag{3.13}. The first of 
them is $J_0$. According to Lemma~\mythelemma{3.3} this Bagderina's 
quantity is related to the quantity $F$ from \mythetag{4.7} introduced
in \mycite{7} by means of the formula \mythetag{3.15}. In our special
coordinates the formula \mythetag{3.15} is immediate from \mythetag{3.14}
due to \mythetag{5.30} and \mythetag{5.31}. Then, using Bagderina's 
formula \mythetag{3.17}, we reproduce \mythetag{3.18}.\par
    The formula \mythetag{3.23} is inapplicable in our special coordinates 
since $\beta^{\,\sssize\text{Bgd}}_1=0$ (see \mythetag{5.30}). There is
a reservation for this case in Bagderina's paper (see Remark 1 on page 7
of \mycite{5}). According to this reservation, if 
$\beta^{\,\sssize\text{Bgd}}_1=0$, then the second formula in 
\mythetag{3.22} is replaced by the following formula:
$$
\hskip -2em
\Cal D^{\,\sssize\text{Bgd}}_2=(J^{\,\sssize\text{Bgd}}_0)^{-4/5}
(\Gamma^{\,\sssize\text{Bgd}}_1\,D_x-\Gamma^{\,\sssize\text{Bgd}}_0
\,D_y).
\mytag{5.34}
$$
Applying the first formula \mythetag{3.22} and the formula \mythetag{5.34},
from \mythetag{3.15}, \mythetag{3.18}, \mythetag{5.30}, and \mythetag{5.33}
we get the following formulas for $\Cal D^{\,\sssize\text{Bgd}}_1$ and
$\Cal D^{\,\sssize\text{Bgd}}_2$:
$$
\xalignat 2
&\hskip -2em
\Cal D^{\,\sssize\text{Bgd}}_1=\frac{B}{F^2}\,D_x,
&&\Cal D^{\,\sssize\text{Bgd}}_2=\frac{3\,F}{B}\,D_y.
\mytag{5.35}
\endxalignat
$$   
The formulas \mythetag{5.35} are in agreement with \mythetag{3.24} and 
\mythetag{3.27} due to \mythetag{3.8}, \mythetag{3.21}, \mythetag{5.11}, 
\mythetag{5.12} and since $D_x$ and $D_y$ are just partial derivatives in 
$x$ and $y$. The formulas \mythetag{5.35} confirm Theorem~\mythetheorem{3.3}. 
\par
     Let's return to the formulas \mythetag{3.13}. As we noted above 
the first of them reduces to \mythetag{3.18}. The others are used in order
to produce Bagderina's scalar invariants \mythetag{4.13}. 
Applying \mythetag{3.13} to \mythetag{4.13}, one can derive some huge
formulas for the invariants $I^{\,\sssize\text{Bgd}}_1$, 
$I^{\,\sssize\text{Bgd}}_2$, $I^{\,\sssize\text{Bgd}}_3$, and
$I^{\,\sssize\text{Bgd}}_4$. In our special coordinates these huge
formulas reduce to the following formulas, which are rather simple:
$$
\gather
\hskip -2em
I^{\,\sssize\text{Bgd}}_1=\frac{15\,F_{\sssize 0.1}}{B}
-\frac{42\,F\,B_{\sssize 0.1}}{5\,B^2}
-\frac{27\,F\,R}{5\,B},
\mytag{5.36}\\
\vspace{1ex}
\hskip -2em
I^{\,\sssize\text{Bgd}}_2=-\frac{6\,B_{\sssize 1.1}}{F\,B}
+\frac{6\,B_{\sssize 1.0}\,B_{\sssize 0.1}}{F\,B^2}
+\frac{9\,R_{\sssize 1.0}}{F}
-\frac{9\,F^4\,S}{B^3},
\mytag{5.37}\\
\vspace{1ex}
\hskip -2em
I^{\,\sssize\text{Bgd}}_3=-\frac{5\,B_{\sssize 1.0}}{F^2}
+\frac{5\,B\,F_{\sssize 1.0}}{F^3},
\mytag{5.38}\\
\vspace{1ex}
\hskip -2em
I^{\,\sssize\text{Bgd}}_4=\frac{90\,F_{\sssize 0.1}}{B}
-\frac{261\,F\,B_{\sssize 0.1}}{5\,B^2}
-\frac{216\,F\,R}{5\,B}.
\mytag{5.39}
\endgather
$$
\mytheorem{5.2} In the case of general position Bagderina's scalar 
invariants $I^{\,\sssize\text{Bgd}}_1$, $I^{\,\sssize\text{Bgd}}_2$, 
$I^{\,\sssize\text{Bgd}}_3$, $I^{\,\sssize\text{Bgd}}_4$ from \mycite{5}
given by the formulas \mythetag{4.13} are expressed through the invariants 
$I_3$, $I_6$, $I_7$, $I_8$ from \mycite{7} by means of the formulas
$$
\gather
\hskip -2em
I^{\,\sssize\text{Bgd}}_3=15\,I_6,\kern 3em
I^{\,\sssize\text{Bgd}}_1=I_3+\frac{14}{5}\,I_8,\kern 3em
I^{\,\sssize\text{Bgd}}_4=3\,I_3+\frac{87}{5}\,I_8,
\quad
\mytag{5.40}\\
\vspace{1ex}
\hskip -2em
I^{\,\sssize\text{Bgd}}_2=-I_7-3\,\bold XI_8+15\,\bold YI_6
+(24\,I_8+15\,I_3)\,I_6, 
\mytag{5.41}
\endgather
$$   
where $\bold X=\Cal D^{\,\sssize\text{Bgd}}_1$ and $\bold Y
=\Cal D^{\,\sssize\text{Bgd}}_2$ are the operators of invariant 
differentiations given by the formulas \mythetag{3.24} and 
\mythetag{3.27}. 
\endproclaim
    The formulas \mythetag{5.40} and \mythetag{5.41} are proved
by direct calculations using the explicit formulas for invariants
\mythetag{5.36}, \mythetag{5.37}, \mythetag{5.38}, \mythetag{5.39}, 
and \mythetag{5.17}. In our special coordinates the operators 
$\Cal D^{\,\sssize\text{Bgd}}_1$ and $\Cal D^{\,\sssize\text{Bgd}}_2$
are given by the formulas \mythetag{5.35}. In the form of vector fields
these operators were introduced in \mycite{7} (see 
Theorem~\mythetheorem{3.3}).\par
     On page 3 of her paper \mycite{5} Yu\.~Yu\.~Bagderina presents
her own comparison formulas expressing $I_3$, $I_6$, $I_7$, $I_8$
through her invariants. We write them as
$$
\gather
\hskip -2em
I_6=\frac{I^{\,\sssize\text{Bgd}}_3}{15},\kern 3em
I_3=\frac{29\,I^{\,\sssize\text{Bgd}}_1}{15}
-\frac{14\,I^{\,\sssize\text{Bgd}}_4}{45},\kern 3em
I_8=\frac{I^{\,\sssize\text{Bgd}}_4}{9}
-\frac{I^{\,\sssize\text{Bgd}}_1}{3},
\quad
\mytag{5.42}\\
\vspace{1ex}
\hskip -2em
I_7=3\,I^{\,\sssize\text{Bgd}}_2-\frac{I^{\,\sssize\text{Bgd}}_5}{3}
\text{, \ where \ }I^{\,\sssize\text{Bgd}}_5
=\frac{J^{\,\sssize\text{Bgd}}_5}{(J^{\,\sssize\text{Bgd}}_0)^{11/5}}.
\mytag{5.43}
\endgather
$$
It is easy to see that the formulas \mythetag{5.42} are converse to the
formulas \mythetag{5.40}. They are easily derived from \mythetag{5.40}.
The formula \mythetag{5.43} is different. It refers the reader to some
quantity $J^{\,\sssize\text{Bgd}}_5$ outside the main Theorem 2 and involves 
a separate bunch of huge notations and calculations on page 25 of
\mycite{5}. Using the relationships \mythetag{5.41} and \mythetag{5.42}, 
we can replace \mythetag{5.43} with the formula
$$
\gathered
I_7=-I^{\,\sssize\text{Bgd}}_2-\frac{1}{3}
\,\Cal D^{\,\sssize\text{Bgd}}_1(I^{\,\sssize\text{Bgd}}_4)
+\Cal D^{\,\sssize\text{Bgd}}_1(I^{\,\sssize\text{Bgd}}_1)\,+\\
+\,\Cal D^{\,\sssize\text{Bgd}}_2(I^{\,\sssize\text{Bgd}}_3)
+\frac{1}{15}\,(21\,I^{\,\sssize\text{Bgd}}_1
-2\,I^{\,\sssize\text{Bgd}}_4)\,I^{\,\sssize\text{Bgd}}_3. 
\endgathered
\mytag{5.44}
$$
\mytheorem{5.3} In the case of general position the scalar invariants
$I_3$, $I_6$, $I_7$, $I_8$ from \mycite{7} given by the formulas
\mythetag{4.4}, \mythetag{4.9}, \mythetag{4.11}, \mythetag{4.12} 
are expressed through Bagderina's scalar invariants $I^{\,\sssize\text{Bgd}}_1$, $I^{\,\sssize\text{Bgd}}_2$, $I^{\,\sssize\text{Bgd}}_3$, $I^{\,\sssize\text{Bgd}}_4$ 
from \mycite{5} by means of the formulas \mythetag{5.42} and \mythetag{5.44},
where $\Cal D^{\,\sssize\text{Bgd}}_1=\bold X$ and \
$\Cal D^{\,\sssize\text{Bgd}}_2=\bold Y$ are the operators of invariant 
differentiations given by the formulas \mythetag{3.24} and \mythetag{3.27}. 
\endproclaim
    Due to Theorems~\mythetheorem{5.2} and \mythetheorem{5.3} the quadruples
of invariants $I_3$, $I_6$, $I_7$, $I_8$ from \mycite{7} and 
$I^{\,\sssize\text{Bgd}}_1$, $I^{\,\sssize\text{Bgd}}_2$, 
$I^{\,\sssize\text{Bgd}}_3$, $I^{\,\sssize\text{Bgd}}_4$ from \mycite{5}
are equivalent to each other.\par
    Apart from $I^{\,\sssize\text{Bgd}}_1$, $I^{\,\sssize\text{Bgd}}_2$, 
$I^{\,\sssize\text{Bgd}}_3$, $I^{\,\sssize\text{Bgd}}_4$, in \mycite{5}
Yu\.~Yu\.~Bagderina considers two omega invariants 
$\Omega^{\,\sssize\text{Bgd}}_1$ and $\Omega^{\,\sssize\text{Bgd}}_2$ 
and writes the relationship 
$$
\hskip -2em
[\Cal D^{\,\sssize\text{Bgd}}_1,\Cal D^{\,\sssize\text{Bgd}}_2]
=\Omega^{\,\sssize\text{Bgd}}_1\,\Cal D^{\,\sssize\text{Bgd}}_1
+\Omega^{\,\sssize\text{Bgd}}_2\,\Cal D^{\,\sssize\text{Bgd}}_2.
\mytag{5.45}
$$ 
Comparing \mythetag{5.45} with \mythetag{5.3} and applying 
Theorem~\mythetheorem{3.3}, we find that 
$$
\xalignat 2
&\hskip -2em
\Omega^{\,\sssize\text{Bgd}}_1=L, 
&&\Omega^{\,\sssize\text{Bgd}}_2=-K.
\mytag{5.46}
\endxalignat
$$
\vskip -1ex
\mytheorem{5.4} In the case of general position Bagderina's omega
invariants $\Omega^{\,\sssize\text{Bgd}}_1$ and $\Omega^{\,\sssize\text{Bgd}}_2$
in \mythetag{5.46} are expressed through the invariants $L$ and $K$ introduced 
in \mycite{7} by means of the formulas \mythetag{5.46}. 
\endproclaim
     For $\Omega^{\,\sssize\text{Bgd}}_1$ and $\Omega^{\,\sssize\text{Bgd}}_2$
from \mythetag{5.46} in item 1 of her Theorem 2 in \mycite{5} Yu\.~Yu\.~Bagderina 
writes the following formulas:
$$
\xalignat 2
&\hskip -2em
\Omega^{\,\sssize\text{Bgd}}_1=\frac{1}{5}(8\,I^{\,\sssize\text{Bgd}}_1
-I^{\,\sssize\text{Bgd}}_4),
&&\Omega^{\,\sssize\text{Bgd}}_2=\frac{1}{5}\,I^{\,\sssize\text{Bgd}}_3.
\mytag{5.47}
\endxalignat
$$
Using \mythetag{5.46} and \mythetag{5.40}, one easily finds that the formulas
\mythetag{5.47} are equivalent to the formulas \mythetag{4.3}.\par
     Apart from the formulas discussed above, which are not new in \mycite{5}
as compared to \mycite{7}, \mycite{3}, and \mycite{4}, in item 1 of her
Theorem 2 Yu\.~Yu\.~Bagderina provides a series of algebraic
equations relating higher order invariants with each other (see (2.9) in
\mycite{5}). These equations are new. But they do not affect the classification
itself. We shall not verify these equations, just saying that this could be 
done with the use of the explicit formulas \mythetag{5.36}, \mythetag{5.37}, 
\mythetag{5.38}, \mythetag{5.39}, and \mythetag{5.35}, which are rather simple. 
In addition to them the equations \mythetag{5.25} and \mythetag{5.26} along
with their compatibility condition \mythetag{5.28} should be taken into 
account.\par
\head
6. Conclusions.
\endhead
     The classification of the equations \mythetag{1.1} suggested by 
Yu\.~Yu\.~Bagderina in \mycite{5} is not absolutely new. It coincides with 
the prior classification suggested in \mycite{3} at least in two items. The 
item 9 of Bagderina's classification Theorem 2 in \mycite{5} coincides with 
the case of {\bf maximal degeneration} from \mycite{3}. The item 1 in 
Bagderina's classification Theorem 2 coincides with the case of {\bf general 
position} from \mycite{3}. The comparison of the two classifications in other 
7 cases (items), which are present in both classifications, will be continued 
in separate papers.\par
     The detailed analysis carried out in the above sections shows that 
in the case of {\bf general position} most structures and most formulas from
Bagderina's  paper \mycite{5} do coincide or are very closely related to those 
in \mycite{7}, though they are given in different notations (see 
Lemma~\mythelemma{3.1}, Lemma~\mythelemma{3.2}, Lemma~\mythelemma{3.3}, 
Theorem~\mythetheorem{3.2}, Lemma~\mythelemma{3.4}, Lemma~\mythelemma{3.5}, 
Theorem~\mythetheorem{3.3}, Theorem~\mythetheorem{5.4}). Four basic scalar 
invariants from item 1 of Bagderina's Theorem 2 in \mycite{5} are equivalent 
to four invariants $I_3$, $I_6$, $I_7$, $I_8$ defined in \mycite{7} (see 
Theorem~\mythetheorem{5.2} and Theorem~\mythetheorem{5.3}).\par
     In addition to the coincidence facts revealed in lemmas and theorems 
listed just above, in section 5 of the present paper we show that most of
the complicated formulas from \mycite{5} do simplify substantially in the
special coordinates suggested in \mycite{7} more than 15 years before
Bagderina's publication.\par

\enddocument
\end